\documentclass{amsart}
\usepackage{amssymb}
\newtheorem{theorem}{Theorem}[section]

\newtheorem{corollary}[theorem]{Corollary}

\newtheorem{example}[theorem]{Example}
\theoremstyle{definition}

\theoremstyle{remark}

\newtheorem{claim}[theorem]{Claim}
\newtheorem{conclusion}[theorem]{Conclusion}

\def\htt{{\rm ht}}

\def\Spec{{\rm Spec}}

\def\Int{{\rm Int}}

\def\Kr{{\rm Kr}}

\begin{document}

\title{Essential Domains and Two Conjectures in Dimension Theory}

\author[M. Fontana]{M. Fontana}

\address{Dipartimento di Matematica, Universit\`a degli Studi ``Roma Tre'', 00146 Roma, Italy}

\email{fontana@mat.uniroma3.it}

\author[S. Kabbaj]{S. Kabbaj}

\address{Department of Mathematics, P.O. Box 5046, KFUPM, Dhahran 31261, Saudi Arabia}

\email{kabbaj@kfupm.edu.sa}

\subjclass[2000]{Primary 13C15, 13F20, 13F05, 13G05, 13B02, 13B30}

\thanks{The first author was partially supported by a research grant
  MIUR 2001/2002 (Cofin 2000-MM01192794). The second author was supported by the Arab Fund for Economic
  and Social Development}

\thanks{This work was done while both authors were visiting Harvard University}

\keywords{Krull dimension, valuative dimension, Jaffard domain, integer-valued polynomial
ring, essential domain, Krull domain, UFD, PVMD, Kronecker function ring, star operation}

\begin{abstract} This note investigates two long-standing conjectures on the Krull dimension
of integer-valued polynomial rings and  of polynomial rings, respectively, in the context of
(locally) essential domains.
\end{abstract}

\maketitle

\section{Introduction}

Let $R$ be an integral domain with quotient field $K$ and let
$\Int(R):=\{f \in K[X]: f(R) \subseteq R\}$ be the ring of
integer-valued polynomials over $R$. Considerable work, part of it
summarized in Cahen-Chabert's book \cite{CC2}, has been concerned
with various aspects of integer-valued polynomial rings. A central
question concerning $\Int(R)$ is to describe its prime spectrum
and, hence, to evaluate its Krull dimension.  Several authors
tackled this problem and many satisfactory results were obtained
for various classes of rings such as Dedekind domains
\cite{Ch1,Ch2}, Noetherian domains \cite{Ch2}, valuation and
pseudo-valuation domains \cite{CH}, and pseudo-valuation domains
of type $n$ \cite{T}. A well-known feature is that $\dim(R[X])-1
\leq \dim(\Int(R))$ for any integral domain $R$ \cite{Ca1}.
However, the problem of improving the upper bound
$\dim(\Int(R))\leq \dim_{v}(R[X])=\dim_{v}(R)+1$ \cite{FIKT1},
where $\dim_{v}(R)$ denotes the valuative dimension of $R$, is
still elusively open in general. It is due, in part, to the fact
that  the fiber in $\Int(R)$ of a maximal ideal of $R$ may have
any dimension \cite[Example 4.3]{Ca1} (this stands for the main
difference between polynomial rings and integer-valued polynomial
rings). Noteworthy is that all examples conceived in the
literature for testing $\dim(\Int(R))$ satisfy the inequality
$\dim(\Int(R)) \leq \dim(R[X])$. In \cite{FIKT1,FIKT2}, we
undertook an extensive study,  under two different approaches, in
order to grasp this phenomenon. We got then further evidence for
the validity of the conjecture:
\begin{center}
 ($\mathcal{C}_{1}$)\quad {\em $\dim(\Int(R))\leq \dim(R[X])$ for any integral domain $R$}.
\end{center}
The current situation can be described as follows; ($\mathcal{C}_{1}$) turned out to be
 true in three large -presumably different- classes of commutative rings; namely: (a)
 Krull-type domains, e.g., unique factorization domains (UFDs) or Krull domains \cite{Gr2,FIKT2};
 (b) pseudo-valuation domains of type $n$ \cite{FIKT2}; and  (c) Jaffard domains \cite{FIKT1,Ca1}.

A finite-dimensional domain $R$ is said to be Jaffard if
$\dim(R[X_{1}, ..., X_{n}]) = n + \dim(R)$ for all $n \geq 1$;
equivalently, if $\dim(R) = \dim_{v}(R)$ \cite{ABDFK,BK,DFK,G,J}. The
class of Jaffard domains contains most of the well-known classes of
finite-dimensional rings involved in dimension theory of commutative rings such as
 Noetherian domains \cite{K}, Pr\"ufer domains \cite{G}, universally catenarian domains
 \cite{BDF}, stably strong S-domains \cite{Kab,MM}. However, the question of establishing
 or denying a possible connection to the family of Krull-like domains is still unsolved. In
 this vein, Bouvier's conjecture (initially, announced during a 1985 graduate course at the
 University of Lyon I) sustains  that:
\begin{center}
($\mathcal{C}_{2}$)\quad {\em finite-dimensional Krull domains, or
more particularly UFDs, need not be Jaffard domains}.
\end{center}
 As the Krull property is stable under adjunction of indeterminates,
 the problem merely deflates to the existence of a Krull domain $R$
 with \(1+\dim(R)\lneqq \dim(R[X])\). It is notable that the rare
 non-Noetherian finite-dimensional UFDs or Krull domains existing in
 the literature do  defeat ($\mathcal{C}_{2}$), since all are Jaffard
 \cite{AM,Da1,Da2,DFK,Fu,G2}. So do the examples of non-Pr\"ufer
 finite-dimensional Pr\"ufer $v$-multiplication domains (PVMDs)
 \cite{Gr1,HMM,MZ,Z}; as a matter of fact, these, mainly, arise as
 polynomial rings over Pr\"ufer domains or as pullbacks, and both
 settings either yield Jaffard domains or turn out to be inconclusive
 (in terms of allowing the construction of counterexamples) \cite{ABDFK,FG}. In order to
  find the missing link, one has then to dig beyond the context of PVMDs.

Essential domains happen to offer such a suitable context for ($\mathcal{C}_{2}$) as
well as a common environment for both conjectures ($\mathcal{C}_{1}$) and ($\mathcal{C}_{2}$),
though these have developed in two dissimilar milieus. An integral domain $R$ is said to be
 {\em essential} if $R$ is an intersection of valuation rings that are localizations of
  $R$ \cite{Gr3}. As this notion does not carry up to localizations, $R$
  is said to be
  {\em locally essential} if $R_{p}$ is essential for each $p\in
  \Spec(R)$. Notice that the locally essential domains correspond
  to the $P$-domains in the sense of Mott and Zafrullah \cite{MZ}. PVMDs and almost Krull
  domains
  \cite[p. 538]{G} are perhaps the
  most important examples of locally essential domains. Recall that Heinzer
  constructed in \cite{H2} an example  of an essential domain
  that is not locally essential. Also, it is worth noticing  that Heinzer-Ohm's example
  \cite{HO} of an essential domain which is not a PVMD is, in fact, locally
  essential (cf. \cite[Example 2.1]{MZ}). Finally recall that a Krull-type
  domain is a PVMD in which no non-zero
  element belongs to an infinite number of maximal $t$-ideals
  \cite{Gr1}. We have thus the following implications within the family of
Krull-like domains:\medskip

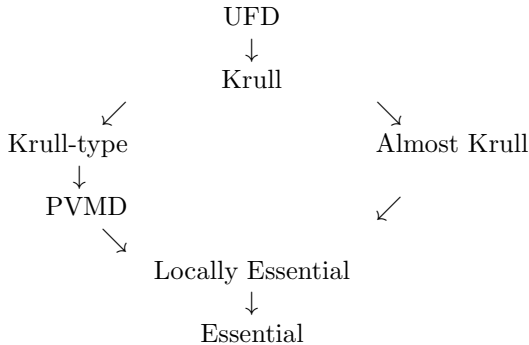
\begin{center}
\begin{tabular}{rcl}
            &UFD\\
            &$\downarrow$\\
            &Krull\\
$\swarrow$  &               &$\searrow$\\
Krull-type  &               &Almost Krull\\
\hfill\hfill $\downarrow$ \hfill\hfill\\
PVMD        &               &$\swarrow$\\
$\searrow$\\
            &Locally Essential\\
            &$\downarrow$\\
            & Essential

\end{tabular}
\end{center}

The purpose of this  note is twofold. First, we state a result
that widens the domain of validity of ($\mathcal{C}_{1}$) to the
class of locally essential domains. It is well-known that
($\mathcal{C}_{1}$) holds for Jaffard domains too
\cite{FIKT1,Ca1}. So one may enlarge the scope of study of
($\mathcal{C}_{2}$) -discussed above- and legitimately raise the
following problem:
\begin{center}
($\mathcal{C'}_2$)\quad {\em Is any finite-dimensional (locally)
 essential domain  Jaffard?}
\end{center}
Clearly, an affirmative answer to ($\mathcal{C'}_2$) will
definitely defeat ($\mathcal{C}_2$); while a negative answer will
partially resolve ($\mathcal{C}_2$) for the class of (locally)
essential domains. Our second aim is to show that the rare
constructions of non-trivial (locally) essential domains (i.e.,
non-PVMD) existing in the literature yield Jaffard domains,
 putting therefore ($\mathcal{C'}_2$) under the status of open problem.
Consequently, a settlement of ($\mathcal{C}_2$) seems -at present-
out of reach.
\section{Result and example}

In the first part of this section, we establish the following result.

\begin{theorem}
For any locally essential domain $R$, $\dim(\Int(R))=\dim(R[X])$.
\end{theorem}

\begin{proof} Assume that $R$ is finite-dimensional and $R\not=K$,
  where $K$ denotes the quotient field of $R$.
Let \(R=\bigcap_{p\in\Delta}R_{p}\) be a locally essential domain,
where $\Delta\subseteq\Spec(R)$. Set:

\(\begin{array}{ll}
\Delta_{1}:= &\{p\in\Delta: R_{p}\ \textup{is a DVR}\}\\
\Delta_{2}:= &\{p\in\Delta: R_{p}\ \textup{is a valuation domain but
                                                          not a DVR}\}.
\end{array}\)\\
We wish to show first that $\dim(\Int(R))\leq\dim(R[X])$. Let $M$ be a
 maximal ideal of $\Int(R)$ such that $\dim(\Int(R))=\htt(M)$ and let
 $\mathcal{M}:= M\cap R$. Without loss of generality, we may assume
 that $\mathcal{M}$ is maximal in $R$ with a finite residue field.  We
 always have $R_{\mathcal{M}}[X]\subseteq (\Int(R))_{\mathcal{M}}\subseteq
 \Int(R_{\mathcal{M}})$ \cite[Corollaires (4), p. 303]{CC1}. If
 $\mathcal{M}\in\Delta_{1}$, then $R_{\mathcal{M}}[X]$ is a two-dimensional
 Jaffard domain \cite[Theorem 4]{S2} and \cite[Proposition 1.2]{ABDFK}. So the
 inclusion $R_{\mathcal{M}}[X]\subseteq(\Int(R))_{\mathcal{M}}$ yields
 $\dim((\Int(R))_{\mathcal{M}})\leq \dim_{v}(R_{\mathcal{M}}[X])=\dim(R_{\mathcal{M}}[X])$.
 Thus, $\dim((\Int(R))_{\mathcal{M}})=\dim(R_{\mathcal{M}}[X])=2$. If
 $\mathcal{M}\in\Delta_{2}$, then $\Int(R_{\mathcal{M}})=R_{\mathcal{M}}[X]=
 (\Int(R))_{\mathcal{M}}$ \cite[Exemples (5), p. 302]{CC1}. If $\mathcal{M}\notin\Delta$,
 then \(R_{\mathcal{M}}=\bigcap_{p\in\Delta, p\subsetneqq\mathcal{M}}R_{p}\) since $R$ is
 a locally essential domain. So that by \cite[Corollaires (3), p. 303]{CC1}
 $\Int(R_{\mathcal{M}})=\bigcap_{p\in\Delta, p\subsetneqq\mathcal{M}}\Int(R_{p})
 =\bigcap_{p\in\Delta, p\subsetneqq\mathcal{M}}R_{p}[X]=R_{\mathcal{M}}[X]=
 (\Int(R))_{\mathcal{M}}$. In all cases, we have $\dim(\Int(R))=
 \dim((\Int(R))_{\mathcal{M}})=\dim(R_{\mathcal{M}}[X])\leq\dim(R[X])$, as desired.

We now establish the inverse inequality
 $\dim(R[X])\leq\dim(\Int(R))$. Let $M$ be a maximal ideal of $R[X]$ such that $\dim(R[X])=\htt(M)$, and
 $\mathcal{M}:= M\cap R$. Necessarily,  $\mathcal{M}$ is maximal in
 $R$. Further, we may assume that $\mathcal{M}$ has a finite residue
 field. If $\mathcal{M}\in\Delta_{2}$ or $\mathcal{M}\notin\Delta$,
 similar arguments as above lead to
 $R_{\mathcal{M}}[X]=(\Int(R))_{\mathcal{M}}$ and hence to the
 conclusion. Next, suppose $\mathcal{M}\in\Delta_{1}$. Then,
 $\Int(R_{\mathcal{M}})$ is a two-dimensional (Pr\"ufer) domain \cite{Ch3}.
 Let $(0)\subsetneqq P_{1}\subsetneqq P_{2}$  be a maximal chain of prime ideals
 in $\Spec(\Int(R_{\mathcal{M}}))$. Clearly, it contracts to
 $(0)\subsetneqq \mathcal{M}R_{\mathcal{M}}$ in $\Spec(R_{\mathcal{M}})$.
 Further, by \cite[Corollaire 5.4]{Ca1}, $P_{1}$ contracts to $(0)$. Therefore,
 by \cite[Proposition 1.3]{Ch2}, $P_{1}=fK[X]\cap\Int(R_{\mathcal{M}})$, for some
 irreducible polynomial $f\in K[X]$. This yields in $\Spec(R_{\mathcal{M}}[X])$ the maximal chain:
\[(0)\subsetneqq P_{1}\cap R_{\mathcal{M}}[X]=fK[X]\cap
R_{\mathcal{M}}[X]\subsetneqq P_{2}\cap R_{\mathcal{M}}[X],\]
which induces in $\Spec(\Int(R))_{\mathcal{M}})$ the following maximal chain:
\[(0)\subsetneqq P_{1}\cap \Int(R))_{\mathcal{M}}=
fK[X]\cap \Int(R))_{\mathcal{M}}\subsetneqq P_{2}\cap \Int(R))_{\mathcal{M}}.\]
Consequently, in all cases, we obtain
$\dim(R[X])=\dim(R_{\mathcal{M}}[X])=\dim((\Int(R))_{\mathcal{M}})$ $\leq\dim(\Int(R))$,
to complete the proof of the theorem.
\end{proof}

>From \cite[Proposition 1.8]{HO} and \cite[Exercise 11, p. 539]{G} we
obtain the following.

\begin{corollary}
Let $R$ be a PVMD or an almost Krull domain. Then $\dim(\Int(R))=\dim(R[X])$.
\end{corollary}

In the second part of this section, we test the problem
($\mathcal{C'}_{2}$) -set up and discussed in the introduction-
for the class of non-PVMD (locally) essential domains. These occur
exclusively in Heinzer-Ohm's example \cite{HO} and Heinzer's
example \cite{H2} both mentioned above. The first of which is a
2-dimensional Jaffard domain \cite[Example 2.1]{MZ}. Heinzer's
example \cite{H2}, too, is a 2-dimensional Jaffard domain by
\cite[Theorem 2.3]{DFK}. Our next example shows that an
enlargement of the scope of this construction -still- generates a
large family of essential domains with nonessential localizations
of arbitrary dimensions $\geq 2$ -but unfortunately- that are
Jaffard domains.


\begin{example}\label{3.1} For any integer $r\geq2$, there
 exists an $r$-dimensional essential Jaffard domain $D$ that is not locally essential.
\end{example}

Notice first that the case $r=2$ corresponds to Heinzer's example
mentioned above. In order to increase the dimension, we modify
Heinzer's original setting by considering Kronecker function rings
via the $b$-operation. For the sake of completeness, we give below
the details of this construction.\medskip

Let $R$ be an integral domain, $K$ its quotient field, $n$ a
positive integer (or $n=\infty$), and $X,X_{1}, ..., X_{n}$
indeterminates over $K$. The $b$-operation on $R$ is the a.b. star
operation defined by $I^{b}:=\bigcap\{IW: W\ \textup{is a
valuation overring of}\ R\}$, for every
 fractional ideal $I$ of $R$. Throughout,  we shall use  $\Kr_{K(X)}(R, b)$ to denote the Kronecker
 function ring of $R$ defined in $K(X)$ with respect to the $b$-operation on $R$; and $R(X_{1}, ..., X_{n})$
  to denote the Nagata ring associated to the polynomial ring $R[X_{1}, ..., X_{n}]$, obtained by
  localizing the latter with respect to the multiplicative system consisting of polynomials
   whose coefficients generate $R$.

Let $r$ be an integer $\geq 2$. Let $k_{0}$ be a field and
$\{X_{n}: n\geq 1\}$, $Y$, $\{Z_{1}, ..., Z_{r-1}\}$ be
indeterminates over $k_{0}$. Let $n$ be a positive integer. Set:
\[\begin{array}{lll}
k_{n}:=k_{0}(X_{1}, ..., X_{n})      &; & k:=\bigcup_{n\geq 1}k_{n}\\
F_{n}:=k_{n}(Z_{1}, ..., Z_{r-1})    &; & F:=\bigcup_{n\geq 1}F_{n}=k(Z_{1}, ..., Z_{r-1})\\
K_{n}:=F_{n}(Y)                      &; & K:=\bigcup_{n\geq 1}K_{n}=F(Y)\\
M_{n}:=YF_{n}[Y]_{(Y)}               &; & M:=\bigcup_{n\geq 1}M_{n}=YF[Y]_{(Y)} \\
A_{n}:= k_{n}+M_{n}                  &; & A:=\bigcup_{n\geq 1}A_{n}=k+M\\
V_{n}:=F_{n}[Y]_{(Y)}                &; & V:=\bigcup_{n\geq
1}V_{n}=F[Y]_{(Y)}.
\end{array}\]

Note that, for each $n\geq 1$, $V$ and $V_{n}$ (resp., $A$ and
$A_{n}$) are one-dimensional discrete valuation domains (resp.,
pseudo-valuation domains) and $\dim_{v}(A)=\dim_{v}(A_{n})$ $=r$.
For each $n\geq 1$, set $X'_{n}:=\frac{1+YX_{n}}{Y}$. Clearly,
$K_{n}=K_{n-1}(X_{n})=K_{n-1}(X'_{n})$. Next, we define
inductively two sequences of integral domains
$\big(B_{n}\big)_{n\geq 2}$ and $\big(D_{n}\big)_{n\geq 1}$  as
follows:
\[\begin{array}{lll}
                                        &;& D_{1}:= A_{1}\\
B_{2}:=\Kr_{K_{1}(X'_{2})}(D_{1},b)     &;& D_{2}:=B_{2}\cap A_{2}\\
B_{n}:=\Kr_{K_{n-1}(X'_{n})}(D_{n-1},b) &;& D_{n}:=B_{n}\cap
A_{n},\textup{ for } n\geq 3.
\end{array}\]

For $ n\geq 2$, let ${\mathcal M}_{n} := M_{n} \cap D_{n}\ (
\subset D_{n} = B_{n}\cap A_{n} \subseteq A_{n}$), where $M_{n}$
is the maximal ideal of $A_{n}$.

\begin{claim} \label{3.4}
(1) $B_{n}$ is an $r$-dimensional Bezout domain.\\
(2) \(B_{n}\cap K_{n-1}=D_{n-1}\subseteq A_{n-1}=A_{n}\cap K_{n-1}\). \\
(3) \(D_{n}\cap K_{n-1}=D_{n-1}\).\\
(4) \(D_{n}[X'_{n}]=B_{n}\) and \((D_{n})_{{\mathcal
M}_{n}}=D_{n}[\frac{1}{YX'_{n}}]=A_{n}\),
with $\frac{1}{X'_{n}}$ and $YX'_{n}\in D_{n}$. \\
(5) ${\mathcal M}_{n}$ is a height-one maximal ideal of $D_{n}$
with
${\mathcal M}_{n}\cap K_{n-1}={\mathcal M}_{n-1}$.\\
(6) For each  $q\in \Spec(D_{n})$ with $q\not={\mathcal M}_{n}$
there exists a unique
prime ideal $Q\in \Spec(B_{n})$ contracting to $q$ in $D_{n}$ and \((D_{n})_{q}=(B_{n})_{Q}\).\\
(7) \(B_{n}=\bigcap\{(D_{n})_{q}: q\in \Spec(D_{n})\ \textup{and}\ q\not={\mathcal M}_{n}\}\).\\
(8)  For each  $q'\in \Spec(D_{n-1})$ with $q'\not={\mathcal
M}_{n-1}$ there exists
 a unique prime ideal $q\ (\not={\mathcal M}_{n})\in \Spec(D_{n})$ contracting to
 $q'$ in $D_{n-1}$ such that \((D_{n})_{q}=(D_{n-1})_{q'}(X'_{n})\).
\end{claim}

\begin{proof} Similar arguments as in \cite{H2} lead to (1)-(7).\\
(8) By (7), $B_{n-1}\subseteq(D_{n-1})_{q'}$, and hence
$(D_{n-1})_{q'}$ is a valuation domain in $K_{n-1}$ of dimension
$\leq r$ containing  $D_{n-1}$. Since $B_{n}$ is the Kronecker
function ring of $D_{n-1}$ defined in $K_{n-1}(X'_{n})$ by all
valuation overrings of $D_{n-1}$, then $(D_{n-1})_{q'}$ has a
unique extension in $K_{n-1}(X'_{n})$, which is a valuation
overring of $B_{n}$, that is, $(D_{n-1})_{q'}(X'_{n})$. By (7),
the center $q$ of $(D_{n-1})_{q'}(X'_{n})$ in $D_{n}$ is the
unique prime ideal of $D_{n}$ lying over $q'$ and
\((D_{n})_{q}=(D_{n-1})_{q'}(X'_{n})\).
\end{proof}

Set $D:=\bigcup_{n\geq 1}D_{n}$ and ${\mathcal M}:=\bigcup_{n\geq
2}{\mathcal M}_{n}$.  It is obvious that $D \subseteq
A=\bigcup_{n\geq 1}A_{n}$.

\begin{claim} \label{3.5} $D_{\mathcal M}=A$ and ${\mathcal M}$ is a height-one maximal
ideal in $D$.
\end{claim}
\begin{proof}
It is an easy consequence of Claim \ref{3.4}(4), since
 ${\mathcal M}_{n} ={\mathcal M}\cap D_{n}$, for each $n$.
\end{proof}

Let $q\in \Spec(D)$ with $q\not={\mathcal M}$.  Then, for some
$m\geq 2$, we have in $D_{m}$, $q_{m}:=q\cap D_{m}\not={\mathcal
M}_{m} ={\mathcal M}\cap D_{m}$.  So, by Claim \ref{3.4}(6),
$B_{m}\subseteq(D_{m})_{q_{m}}$, hence $(D_{m})_{q_{m}}$ is a
valuation overring of $B_{m}$ of dimension $\leq r$, whence, by
Claim~\ref{3.4}(8), $D_{q}=\bigcup_{n\geq
1}(D_{n})_{q_{n}}=(D_{m})_{q_{m}}(X'_{m+1}, ...)$ is still a
valuation domain of dimension $\leq r$.

\begin{claim} \label{3.6}
$D=\bigcap\{D_{q}: q\in \Spec(D)\ \textup{and}\ q\not={\mathcal
M}\}$.
\end{claim}

\begin{proof}
Similar to \cite{H2}.
\end{proof}


>From Claims \ref{3.5} and \ref{3.6} we obtain:

\begin{conclusion}
$D$ is an essential domain with a nonessential localization and
$\dim(D)=\dim_{v}(D)=r$.
\end{conclusion}

\end{document}